\documentclass[12pt]{article}
\usepackage{amsthm}
\usepackage{amsmath}
\usepackage{amsfonts}
\usepackage{amssymb}
\usepackage{url}
\usepackage{epsfig}

\newtheorem{Thm}{Theorem}[section]
\newtheorem{Cor}[Thm]{Corollary}
\newtheorem{Lem}[Thm]{Lemma}

\newtheorem{Conj}[Thm]{Conjecture}

\theoremstyle{definition}
\newtheorem{Def}[Thm]{Definition}
\newtheorem{Xmp}[Thm]{Example}
\newtheorem{Rem}[Thm]{Remark}
\newtheorem{Alg}[Thm]{Algorithm}
\newtheorem{Prob}[Thm]{Problem}

   \def\Prob{{\bf P}}

   \def\C{{\mathbb C}}

   \def\R{{\mathbb R}}
   \def\N{{\mathbb N}}
   \def\S{{\mathbb S}}

   \def\S{{\Sigma}}

   \def\Pf{{\it Proof.$\;\;$}}

   \def\qed{\hfill$\diamond$}

   \def\cone{\mbox{\rm cone~}}
   
   \def\det{\mbox{\rm det~}}
   
   \def\diag{\mbox{\rm diag~}}
   \def\image{\mbox{\rm image~}}
   
   \def\rk{\mbox{\rm rk~}}
   \def\SLC{\mbox{\rm SLC~}}
   
   \def\spann{\mbox{\rm span}}
   \def\tr{\mbox{\rm tr~}}

   \def\({\langle}
   \def\){\rangle}
   \def\mb{\boldsymbol}

   \def\cM{{\mathcal M}}

   \def\cP{{\mathcal P}}

   \def\cS{{\mathcal S}}

   \def\1{\mb1}

   \def\v0{{\bf 0}}

\date{}

\begin{document}
\bibliographystyle{plain}
\title{\bf Equations for hidden Markov models}
\author{Alexander Sch\"onhuth}

\maketitle

\begin{abstract}
We will outline novel approaches to derive model invariants for hidden
Markov and related models. These approaches are based on a theoretical
framework that arises from viewing random processes as elements of the
vector space of string functions. Theorems available from that
framework then give rise to novel ideas to obtain model invariants
for hidden Markov and related models.
\end{abstract}

\section{Introduction}
\label{sec.intro}

In the following, we will outline how to obtain invariants for hidden
Markov and related models, based on an approach which, in its most
prevalent application, served to solve the identifiability problem for
hidden Markov processes (HMPs) in 1992 \cite{Ito92}. Some of its
foundations had been layed in the late 50's and early 60's in order to
get a grasp of problems related to that of identifying HMPs
\cite{Blackwell57,Gilbert59,Dharma63a,Dharma63b,Dharma65,Heller65}.
The approach can be viewed as being centered around the definition of
{\em finite-dimensional} discrete-time, discrete-valued stochastic
processes (referred to as {\em discrete random processes} in the
following)\footnote{In the literature, finite-dimensional discrete
  random processes are alternatively referred to as {\em finitary}
  \cite{Heller65} or {\em linearly dependent} \cite{Ito92}
  processes. In the following, we will stay with the term
  finite-dimensional (discrete random processes) in accordance with
  the latest contributions on the topic
  \cite{Jaeger00,Faigle07,Schoenhuth07,Schoenhuth08,Schoenhuth09}}. It
Examples of finite-dimensional discrete random processes other than
HMPs are quantum random walks (QRWs). QRWs have been brought up mostly
to emulate Markov chain related algorithms (e.g.~Markov Chain Monte
Carlo techniques) on quantum computers \cite{Aharonov01}.\par In the
following, we will introduce finite-dimensional string functions and
formally describe how to view discrete random processes as string
functions. We will further provide helpful characterizations and
related theorems. In sec.~\ref{sec.finiteinvariants} we will determine
polynomials that generate the ideal of the invariants of the
finite-dimensional model, in the usual sense of algebraic
statistics. In sec.~\ref{sec.finitedet}, we will prove a theorem from
which, as a corollary, one obtains a proof of conjecture 11.9 in
\cite{Bray05}. This corollary will be listed in sec.~\ref{sec.hmms}
where we will draw the connections to the hidden Markov model in more
detail. In sec.~\ref{sec.markov} we will show how to obtain invariants
for the Markov models, based on the results of the preceding sections.
In sec.~\ref{sec.tracealg} we will briefly demonstrate that trace
algebras, as well, can be viewed as certain finite-dimensional string
functions. Invariants of the finite-dimensional model are relatively
easy to obtain

\section{Preliminaries: String Functions}
\label{sec.prelim}

Detailed proofs and explanations of the following results can be found
from \cite{Schoenhuth07}.  Let $\Sigma^*=\cup_{n\ge 0}\Sigma^n$ denote
the set of all strings of finite length over the finite alphabet
$\Sigma$ where the word $\Box\in\Sigma^0$ of length $|\Box|=0$ is the
{\em empty string}. Single letters are usually denoted by $a,b$
whereas strings of arbitrary length are denoted by $v,w$ (for example,
$v=a_1...a_n\in\S^n,w=b_1...b_m\in\S^m$ where $a_i,b_j\in\S$). We have
the concatenation operation:
\begin{equation}
w\in \Sigma^m, v\in \Sigma^n \quad \Longrightarrow\quad wv\in\Sigma^{m+n}.
\end{equation}
We denote the \emph{length} of $v\in\Sigma^n$ by $|v|=n$. We
now direct our attention to real-valued string functions
\begin{equation}
p:\;\S^*\longrightarrow\R
\end{equation}
and further to $\R^{\S^*}$, that is, to the real vector space of string
functions over $\S$. The notation $p$ is due to that discrete random
processes will be viewed as string functions, which will be described
in the following.

\subsection{Discrete Random Processes as String Functions}
\label{ssec.processes}

\begin{sloppypar}
Given a discrete random process $(X_t)$ with values in the alphabet
$\S$, the prescription
\begin{equation*}
p_X(v=a_1...a_n) = \Prob (\{X_1=a_1,...,X_n=a_n\})
\end{equation*}
gives rise to a string function $p_X$ associated with the random
process.  $p_X(a_1...a_n)$ then just is the probability that the
associated random process emits the string $a_1...a_n$ at periods
$t=1,...,n$. String functions associated with discrete random processes
can be characterized as follows.\end{sloppypar}

\begin{Thm}
\label{t.ssf}
A string function $p:\S^*\to\R$ is associated with a discrete
random process iff the following conditions hold.
\begin{enumerate}
\item[(a)] $p(v)\ge 0$ for all $v\in\S^*$.
\item[(b)] $\sum_{a\in\S}p(va) = p(v)$ for all $v\in\S^*$.
\item[(c)] $p(\Box) = 1$.
\end{enumerate}
\end{Thm}

Note that $(b)$ in combination with $(c)$ implies 
\begin{equation}
\forall n\ge 0:\quad\sum_{v\in\S^n}p(v) = 1.
\end{equation}

\begin{Def}
\label{d.ssf}
A string function $p:\S^*\to\R$ is called 
\begin{itemize}
\item {\em stochastic string
  function (SSF)} if it is associated with a discrete random process, that is,
  iff $(a),(b)$ and
  $(c)$ of theorem \ref{t.ssf} apply,
\item {\em unconstrained stochastic string function (USSF)} if only
  $(a)$ and $(b)$ apply (in accordance with the terminology of
  \cite{Pachter05}) and
\item {\em generalized unconstrained stochastic string function
  (GUSSF)} if only $(b)$ applies. 
\end{itemize}
\end{Def}

In the following, the terms (generalized unconstrained) random process
and (GU)SSF will be used interchangeably. Furthermore, note that
$p(a_1...a_n)$ just is a different notation for $p_{a_1...a_l}$ which
was used in \cite{Pachter05}. \par

\subsection{Dimension of String Functions}
\label{ssec.dim}

The following definitions are fundamental for this work.

\begin{Def}
\label{d.hankel}
Let $p:\S^*\to\R$ be a string function over $\S$. Then
\begin{equation}
\cP_p := [p(wv)_{v,w\in \Sigma^*}]\in\R^{\Sigma^*\times\Sigma^*}
\end{equation}
is called the {\em Hankel matrix} of $p$ (also called {\em prediction
  matrix} in case of a SSF $p$). We define 
\begin{equation}
\dim p := \rk\cP_p
\end{equation}
to be the {\em dimension} of $p$. In case of $\dim p<\infty$
the string function $p$ is said to be
{\em finite-dimensional}.
\end{Def}

\begin{Xmp} Let $p:\S^*\to\R$ be a string function over the binary alphabet $\S=\{0,1\}$.
\begin{equation*}
\cP_p = 
  \begin{pmatrix} 
       p(\Box) & p(0) & p(1)   & p(00)   & p(01)   & p(10)   & p(11)   & \hdots \\
       p(0)  & p(00)  & p(10)  & p(000)  & p(010)  & p(100)  & p(110)  & \hdots \\
       p(1)  & p(01)  & p(11)  & p(001)  & p(011)  & p(101)  & p(111)  & \hdots \\
       p(00) & p(000) & p(100) & p(0000) & p(0100) & p(1000) & p(1100) & \hdots \\
       p(01) & p(001) & p(101) & p(0001) & p(0101) & p(1001) & p(1101) & \hdots \\
       p(10) & p(010) & p(110) & p(0010) & p(0110) & p(1010) & p(1110) & \hdots \\
       p(11) & p(011) & p(111) & p(0011) & p(0111) & p(1011) & p(1111) & \hdots \\
       \vdots& \vdots & \vdots & \vdots  & \vdots  & \vdots  & \vdots  & \ddots
  \end{pmatrix}
\end{equation*}
then is the Hankel matrix where strings of finite length have been
ordered lexicographically. Note that within a row values refer to
strings that have the same suffix whereas within a column values refer
to strings that have the same prefix. See also \cite{Finesso06} for an
example of a Hankel matrix.
\end{Xmp} 

The following characterization of finite-dimensional string functions is the major
source of motivation for this work.

\begin{Thm}[\cite{Ito92,Jaeger00,Schoenhuth07}]
\label{t.oom}
Let $p:\S^*\to\R$ be a string function.
Then the following conditions are equivalent.
\begin{enumerate}
\item[(i)] $p$ has dimension at most $d$.
\item[(ii)] There exist vectors $x,y\in\R^d$ as well as matrices
            $T_a\in\R^{d\times d}$ for all $a\in\S$ such that
\begin{equation}
\label{eq.oom}
\forall v\in\S^*:\quad p(v=a_1...a_n) = \langle y|T_{a_n}...T_{a_1}|x\rangle.
\end{equation}
\end{enumerate}
\end{Thm}

Fully elaborated proofs of theorem \ref{t.oom} can be found in
\cite{Jaeger00,Schoenhuth07}.  Note that (\ref{eq.oom}) can be
transformed to
\begin{equation}
p(v) = \tr T_{a_n}...T_{a_1}C
\end{equation}
where $C=xy^T\in\R^{d\times d}$.\par 

\begin{Xmp}
\label{x.hmm}
The most prominent example for finite-dimensional SSFs are hidden
Markov chains.  Let $p$ be an SSF associated with a hidden Markov
chain on $d$ hidden states and output alphabet $\S$. Let
$A=(\Prob(i\to j))_{1\le i,j\le d}$ be the transition probability
matrix, $E_{ia},1\le i\le d,a\in\S$ be the emission probabilities and
$\pi$ be the initial probability distribution. We define 
\begin{equation}
O_a:=\diag(E_{ia},i=1,...,d)\in\R^{d\times d}
\end{equation}
and further 
\begin{equation*}
T_a := A^TO_a\in\R^{d\times d}.
\end{equation*}
The $T_a$ together with $y:=(1,...,1)\in\R^d$, $x:=\pi\in\R^d$ then
provide a representation corresponding to (\ref{eq.oom}).  
\end{Xmp}

We will be particularly interested in finite-dimensional GUSSFs (we
recall definition \ref{d.ssf}). The following theorem provides a
characterization.

\begin{Thm}
\label{t.finitegussf}
Let $p:\S^*\to\R$ be a string function such that $\dim p\le d$. Then
the following two statments are equivalent:
\begin{enumerate}
\item[(i)] $p$ is a GUSSF, that is, $\sum_{a\in\S}p(va)=p(v)$ for all $v\in\S^*$.
\item[(ii)] There exist vectors $x,y\in\R^d$ as well as matrices
            $T_a\in\R^{d\times d}$ for all $a\in\S$ such that
\begin{equation}
\label{eq.oom2}
\forall v\in\S^*:\quad p(v=a_1...a_n) = \langle y|T_{a_n}...T_{a_1}|x\rangle.
\end{equation}
as well as
\begin{equation}
\label{eq.gussf}
y^T\sum_{a\in\S}T_a = y^T 
\end{equation}
translating to that $y$ is an eigenvector of the eigenvalue $1$ of the
transpose of $\sum_{a\in\S}T_a$.\\
\end{enumerate}
\end{Thm}

\noindent In the following, we will write
\begin{equation}
\label{eq.Tv}
T_v := T_{a_n}...T_{a_1}, T_w=T_{b_m}...T_{b_1}
\end{equation}
in case of $v=a_1...a_n\in\S^n,w=b_1...b_m\in\S^m$.\\

\noindent\Pf The obvious direction is ``$\Leftarrow$'':
\begin{equation}
\begin{split}
\sum_{a\in\S}p(va) &\stackrel{(\ref{eq.oom2})}{=} \sum_{a\in\S}\langle y|T_aT_v|x \rangle
= \langle y|\sum_{a\in\S}T_a|T_vx \rangle \\
&\stackrel{(\ref{eq.gussf})}{=} \langle y|T_v|x \rangle= p(v).
\end{split}
\end{equation}
For ``$\Rightarrow$'', let $d^*:=\dim p\le d$ be the actual dimension
of $p$.  According to theorem \ref{t.oom}, we find matrices
$\tilde{T_a}\in\R^{d^*\times d^*},a\in\S$ and vectors
$\tilde{x},\tilde{y}\in \R^{d^*}$ such that
\begin{equation}
\label{eq.gussf2}
\forall z\in\R^{d^*}:\quad \langle
\tilde{y}|\sum_{a\in\S}\tilde{T_a}|z \rangle = \langle \tilde{y}|z \rangle.
\end{equation}
In case of $d^*=d$ we will have proven the claim by putting
$T_a:=\tilde{T_a},x=\tilde{x},y=\tilde{y}$. In case of $d^*<d$ we will
obtain suitable matrices $T_a\in\R^{d\times d}$ and vectors $x,y\in\R^d$ by 
putting
\begin{equation}
(T_a)_{ij},x_i,y_i := \begin{cases} (\tilde{T_a})_{ij},\tilde{x}_i,\tilde{y}_i & 1\le i,j\le d^*\\
                                   0 & \text{else}
                     \end{cases}.
\end{equation}
From theorem \ref{t.oom} we obtain matrices
$\tilde{T_a}\in\R^{d^*\times d^*},a\in\S$ and vectors
$\tilde{x},\tilde{y}\in \R^{d^*}$ such that
\begin{equation}
\label{eq.oom3}
p(v=a_1...a_n) = \langle y|T_v|x \rangle.
\end{equation}
Condition $(i)$ then implies that
\begin{equation}
\langle y|(\sum_{a\in\S}T_a)T_v|x \rangle = \sum_{a\in\S}\langle y|T_aT_v|x \rangle 
\stackrel{(\ref{eq.oom3})}{=} \sum_{a\in\S}p(va) \stackrel{(i)}{=} p(v) = \langle y|T_v|x \rangle.
\end{equation}
It remains to show that  
\begin{equation}
\spann\{T_vx\,|\,v\in\S^*\} = \R^{d^*}.
\end{equation}
However,
assuming the contrary would lead to the contradiction
\begin{multline}
d^* = \dim p = \rk [p(wv)]_{v,w\in\S^*} = \rk [\langle y|T_vT_w|x \rangle]_{v,w\in\S^*}\\
 \le \dim\spann\{T_wx\,|\,w\in\S^*\} < d^*.\quad\diamond
\end{multline}

\noindent Matrices $T_a$ can be
computationally determined according to a procedure which we will
describe in the following. Therefore, for a string function $p$,
we introduce the notation
\begin{equation}
\begin{array}{rccc}
 p_v:&\S^*&\to&\R\\
     &w&\mapsto&p(wv)
\end{array}\;\text{resp.}\;
\begin{array}{rccc}
 p^w:&\S^*&\to&\R\\
     &v&\mapsto&p(wv)
\end{array}.
\end{equation}
That is, the $p_v$ resp.~$p^v$ are the row resp.~column vectors of the
Hankel matrix $\cP_p$. These are string functions in their own
right. Note that $p_{\Box}=p^{\Box}=p$. Moreover, note that in case of
a stochastic process $p$ s.t.~$p(w=b_1...b_{m})\ne 0$ it holds that
\begin{multline}
\frac{1}{p(w)}p^w(v=a_1...a_l)\\ 
= \Prob(\{X_{l'+1}=a_1,...,X_{l'+l}=a_l\}\,|\,\{X_1=b_1,...,X_{l'}=b_{l'}\}).
\end{multline}
Therefore, $\frac{1}{p(w)}p^w$ is just the discrete random process
being governed by the probabilities of $p$ conditioned on that $w$ has
already been emitted.\\

The following is a generic algorithmic strategy to infer matrices
$T_a\in\R^{d\times d}$ and $x,y\in\R^d$ corresponding to
(\ref{eq.oom}) from a finite-dimensional Hankel matrix. At this point,
the algorithm needs the entire string function $p$ as an input. We
will explain how to obtain a practical version of this generic
strategy later in this section.

\begin{Alg}
\label{a.ooms}
$     $\\[-.5ex]
\hrule
\medskip
\noindent{\bf Input}: A string function $p$ such that $\dim p = d<\infty$.\\[.5ex]
\noindent{\bf Output}: Matrices $T_a\in\R^{d\times d},a\in\S$ and vectors $x,y\in\R^d$
such that 
\begin{equation}
p(v=a_1...a_n) = \tr T_{a_n}...T_{a_1}xy^T.
\end{equation}
\hrule
\medskip 
\begin{enumerate}
\item Determine words $v_1,...,v_d$ resp.~$w_1,...,w_d$ such that the
  $f_{v_i}$ resp.~the $g_{w_j}$ span the row resp.~column space of
  $\cP_p$. Hence the matrix
\begin{equation}
\label{eq.V}
V := [p(w_jv_i)]_{1\le i,j\le d}
\end{equation}
has full rank $d=\dim p$.
\item Denote by $V_i$ resp.~$V^j$ the $i$-th row resp.~the $j$-th
  column of $V$ and define
\begin{equation}
x=(x_1,...,x_d)^T := (p(v_1),...,p(v_d))^T
\end{equation}
and $y=(y_1,...,y_d)\in\R^d$ such that 
\begin{equation}
(p(v_1),...,p(v_d)) = \sum_{i=1}^dy_iV_i
\end{equation}
which can be done as $p_{\Box}=p$ (the uppermost row of the Hankel
matrix) is linearly dependent of the $p_{v_i}$ (the basis of the 
row space of the Hankel matrix).
\item For each $a\in\S$, determine matrices
\begin{equation}
\label{eq.Wa}
W_a:= [p(w_jav_i)]_{1\le i,j\le d}.
\in\R^{d\times d}.
\end{equation}
\item One can then show that $x,y$ and $T_a:=(W_aV^{-1}),a\in\S$
are as needed for theorem \ref{t.oom}.
\end{enumerate}
\end{Alg}

Clearly, the driving question behind algorithm \ref{a.ooms}
is its practicability. A first clue to this is the following theorem.
Therefore, we set 
\begin{equation}
\S^{\le n}:=\cup_{t=0}^{n}\S^t
\end{equation}
to be the
set of all strings of length at most $n$ and
\begin{equation}
\label{eq.hankelminor}
\cP_{p,n,m}:=[p(wv)]_{|v|\le n,|w|\le m}\in\R^{\S^{\le n}\times\S^{\le m}}.
\end{equation}
to be the finite minor of the Hankel matrix referring to row
resp.~column vectors indexed by strings of length at most $n$
resp.~$m$.

\begin{Thm}
\label{t.dimcheck}
Let $p:\S\to\R$ be a string function such that $\dim p\le d$. Then
it holds that
\begin{equation}
\dim p = \rk\cP_{p,d-1,d-1}.
\end{equation}
\end{Thm}

This means that, given an upper bound $d$ on the dimension of $p$, the
dimension of $p$ can be determined by inspecting the
finite-dimensional matrix $\cP_{p,d-1,d-1}$. See \cite{Schoenhuth07}
for a proof. Note, however, that the size of $\cP_{p,d-1,d-1}$ is
exponential in $d$ such that naive approaches to determining $V$
(\ref{eq.V}) would result in exponential runtime. The final clue to
the practicability of algorithm \ref{a.ooms} is an efficient algorithm
to determine $V$ which has recently been presented
\cite{Schoenhuth08}. The algorithm applies in case one is provided
with an arbitrary generating system of the row or column space of
$\cP_p$. Corresponding generating systems emerge naturally for
finite-dimensional processes of interests, in particular for hidden
Markov processes and also for quantum random walks.\par A consequence
of theorem \ref{t.dimcheck} is

\begin{Thm}[\cite{Schoenhuth07}]
\label{t.unique}
Let $p$ be a string function such that $\dim p\le d$. Then $p$
is uniquely determined by the values
\begin{equation}
p(v), \quad |v|\le 2d-1.
\end{equation}
\end{Thm}

\begin{sloppypar}
{\em Proof Sketch.} The idea is, given two string functions $p_1,p_2$
where $\dim p_1,\dim p_2\le d$ which coincide on strings of length up
to $2d-1$, to determine matrices $T_a$ and vectors $x,y$ as in theorem
\ref{t.oom} according to algorithm \ref{a.ooms} for both $p_1$ and
$p_2$. Thanks to theorem \ref{t.dimcheck}, in algorithm \ref{a.ooms},
$V$ can be determined by inspecting values of $p$ at strings of length
at most $2d-2$ in $\cP_{p,d-1,d-1}$ and, subsequently, by inspecting
strings of length at most $2(d-1)+1=2d-1$ in order to obtain the
$W_a$. As $p_1$ and $p_2$ coincide on strings of length $2d-1$, this
will result in the same $V$ and $W_a$.  Hence $p_1=p_2$.\qed\\
\end{sloppypar}

The following corollary is an obvious consequence of
theorem \ref{t.unique} due to property $(b)$ from
theorem \ref{t.ssf}. However, it had been well-known already before.  See
e.g.~\cite{Merhav02,Finesso06} and the references therein.

\begin{Cor}[\cite{Merhav02,Finesso06}]
\label{c.unique}
A GUSSF $p$ such that $\dim p\le d$ is uniquely determined by the values 
\begin{equation}
p(v), \quad |v| = 2d-1.
\end{equation}
In other words, a discrete random process whose dimension can be upper
bounded by $d$ is uniquely determined by its probability distribution
over the strings of length $2d-1$.
\end{Cor}

\begin{Rem}
Note that for a string function $p$ with $\dim p\le d<\infty$, rows
and columns of the Hankel matrix indexed by strings of length at least
$d$ must necessarily be linearly dependent of their counterparts
referring to strings of length at most $d-1$. These observations are
crucial for the core result of the following section.
\end{Rem}

\section{Finite-Dimensional Models}
\label{sec.finiteinvariants}

Finite-dimensional models over $\S$ are defined to be the polynomial maps
\begin{equation}
\label{eq.gnd}
\begin{array}{rccc}
\mathbf{g}_{n,d}: & \cS^d\subset\C^{|\S|d^2+2d} & \longrightarrow & \C^{|\S|^n}\\
                 & ((T_a)_{a\in\S}),x,y) & \mapsto & (\langle y|T_{a_n}...T_{a_1}|x \rangle)_{v=a_1...a_n\in\S^n}.
\end{array}
\end{equation}
where 
\begin{equation}
((T_a)_{a\in\S}),x,y)\in\cS^d\quad :\Leftrightarrow\quad 
y^T\sum_{a\in\S}T_a = y^T.
\end{equation}
According to theorem \ref{t.finitegussf}, $\cS^d$ comprises precisely
the parameterizations of the generalized unconstrained random
processes of dimension at most $d$. Obviously,
\begin{equation}
\cS^d\cong\C^{(|\S|-1)d^2+d(d-1)+2d}.
\end{equation}
Therefore, the Zariski closure of $\image(\mathbf{g}_{n,d})$ is an irreducible
variety.\\ 

In the following, we will make use of the polynomial map
(\ref{eq.gnd}) to derive a set-theoretic theorem with a strong view
towards the invariants of the Zariski closure of the image of
$\mathbf{g}_{n,d}$.  In case of $n\ge 2d-1$, invariants
for the image of ${\mathbf{g}_{n,d}}$ can be derived by
inspection of the Hankel matrix.  As in (\ref{eq.hankelminor}), let
$\cP_{p,n,m}$ be the partial Hankel matrix that is filled with all
values $p(wv)$ such that $|v|\le n,|w|\le m$.

\begin{Thm}
\label{t.generators}
Let $n\ge 2d-1$ and $(p(v))_{v\in\S^n}$ be an (unconstrained)
probability distribution. Then it holds that
\begin{equation*}
(p(v))_{v\in\S^n}\in\image(\mathbf{g}_{n,d})
\end{equation*}
if and only if the following two conditions apply where, in case of $|u|<n$,
\begin{equation}
\label{eq.expand}
p(u)=\sum_{u\in\S^{n-k}}p(uv).
\end{equation}
\begin{enumerate}
\item[(a)]
\begin{equation}
\label{eq.genaa}
\det[p(w_jv_i)]_{1\le i,j\le d+1} = 0
\end{equation}
for all choices of words $v_1,...v_{d+1},w_1,...w_{d+1}$ of length at most
$d-1$, which can be equivalently put as
\begin{equation}
\label{eq.genab}
\rk\cP_{p,d-1,d-1}\le d
\end{equation}
\item[(b)] 
\begin{equation}
\label{eq.condb}
\rk \cP_{p,\lceil\frac{n}{2}\rceil, \lfloor\frac{n}{2}\rfloor}
= \rk \cP_{p,\lfloor\frac{n}{2}\rfloor, \lceil\frac{n}{2}\rceil }
= \rk \cP_{p,d-1,d-1}
\end{equation}
\end{enumerate}
\end{Thm}

(\ref{eq.condb}) states that rows resp.~columns in
$\cP_{p,\lceil\frac{n}{2}\rceil, \lfloor\frac{n}{2}\rfloor}$ and
$\cP_{p,\lfloor\frac{n}{2}\rfloor, \lceil\frac{n}{2}\rceil }$
referring to row strings $v$ resp.~column strings $w$ where $|v|,|w|
\ge d$ are linearly dependent of their counterparts in
$\cP_{p,\lceil\frac{n}{2}\rceil, \lfloor\frac{n}{2}\rfloor}$ and
$\cP_{p,\lfloor\frac{n}{2}\rfloor, \lceil\frac{n}{2}\rceil }$ that
refer to row resp.~column strings of length at most $d-1$.\\

\Pf ``$\Rightarrow$'': Let $(p(v))_{v\in\S^n}$ be in the image of
$\mathbf{g}_{n,d}$. Theorem \ref{t.oom} states that the Hankel matrix
$\cP_p$ of $p$ has rank at most $d$. This implies $(a)$ as it just
expresses that some Hankel matrix minors of size $d+1$ do not have full
rank.\par  Theorem \ref{t.dimcheck} then states that bases of the row
resp.~the column space of $\cP_p$ can be obtained by inspecting
row resp.~column vectors referring to strings of length at most
$d-1$ which implies $(b)$.\\

``$\Leftarrow$'': Let $(p(u))_{u\in\S^n}$ s.t.~$(a),(b)$ apply.  In
order to prove that $(p(u))_{u\in\S^n}\in\image\mathbf{g}_{n,d}$, we
have to provide a parameterization $((T_a)_{a\in\S},x,y)\in\cS^d$ such
that
\begin{equation}
\label{eq.oomrepr}
p(u=a_1...a_n)=y^TT_{a_n}...T_{a_1}x
\end{equation}
\begin{sloppypar}
\noindent for all strings $u\in\S^n$.  Therefore, we will
provide a parameterization $((T_a)_{a\in\S},x,y)\in\C^{|\S|d^2+2d}$ such
that\end{sloppypar}
\begin{equation}
\label{eq.oomrepr2}
p(u=a_1...a_k)=y^TT_{a_k}...T_{a_1}x
\end{equation}
for all strings $u$ such that $|u|\le n$ where $p(u)$ is defined
according to (\ref{eq.expand}) in case of $|u|<n$.  By this definition
of $p(u),|u|<n$, it is straightforward to show that
$((T_a)_{a\in\S},x,y)\in\cS^{d}$ which completes the proof.
Furthermore, note that it suffices to provide a parameterization
$((T_a)_{a\in\S},x,y)\in\cS^{d^*}$ for arbitrary $d^*\le d$ since, in
case of $d^*<d$, we extend the $T_a$ as well as $x,y$ by zero entries
to obtain a $d$-dimensional parametrization from $\cS^d$.  Combining
these facts, we have to show that, for suitable $d^*\le d$, there are
matrices $T_a\in\R^{d^*\times d^*}$ and vectors $x,y\in\R^{d^*}$ such
that (\ref{eq.oomrepr2}) holds.\par We obtain the desired
parameterization $((T_a)_{a\in\S},x,y)$ according to the ideas of
algorithm \ref{a.ooms}. First, determine strings $v_1,...,v_{d^*}$ and
$w_1,...,w_{d^*}$ of length at most $d-1$ such that
\begin{equation}
V:=[p(w_jv_i)]_{1\le i,j\le d^*}
\end{equation}
has full rank $d^*:=\rk\cP_{p,d-1,d-1}\le d$. We define
\begin{equation}
x=(x_1,...,x_{d^*})^T := (p(v_1),...,p(v_{d^*}))^T
\end{equation}
and $y=(y_1,...,y_{d^*})\in\R^{d^*}$ such that 
\begin{equation}
(p(v_1),...,p(v_{d^*})) = \sum_{i=1}^{d^*}y_iV_i
\end{equation}
where $V_i=(p(v_iw_1),...,p(v_iw_{d^*})^T$ is the $i$-th row of $V$
which can be done since the uppermost row of $\cP_{p,n,d-1}$ is
linearly dependent of the rows referring to the strings $v_i$.
Furthermore, for each $a\in\S$, we determine matrices
\begin{equation}
W_a:= [p(w_jav_i)]_{1\le i,j\le d^*}.
\in\R^{d^*\times d^*}
\end{equation}
Note that probabilities in $W_a$ may refer to strings of length up to
$2d-1$ which establishes the necessity of the assumption $n\ge 2d-1$.
We then claim that defining
\begin{equation}
T_a:=W_aV^{-1}
\end{equation}
gives rise to the desired parametrization in terms of
(\ref{eq.oomrepr2}). We will obtain an easy proof of this claim
by three elementary lemmata.\par

\begin{Lem} 
\label{l.generators.1}
For all $v,w\in\S^*$ such that $|wv|\le\lceil\frac{n}{2}\rceil$
($T_v=T_{a_k}...T_{a_1},v=a_1...a_k\in\S^k$):
\begin{equation}
T_v\begin{pmatrix} p(wv_1)\\ \vdots\\ p(wv_{d^*})
   \end{pmatrix}
= \begin{pmatrix} p(wvv_1)\\ \vdots\\ p(wvv_{d^*})
   \end{pmatrix}
\end{equation}
\end{Lem}

{\bf Proof of lemma \ref{l.generators.1}}: Note first that $|v_i|\le
d-1\le \frac{2d-1}{2}\le \frac{n}{2}$ which implies $|wvv_i|\le n$. As
$(p(wv_1),...,p(wv_{d^*}))^T$ is contained in the column space of $V$
it suffices to show the statement for $w=w_j$.  We do this by
induction on $|v|$:\par
\noindent $|v|=1$: 
\begin{equation}
\begin{split}
T_a\begin{pmatrix} p(w_jv_1)\\ \vdots\\ p(w_jv_{d^*})
   \end{pmatrix} 
&= W_aV^{-1} \begin{pmatrix} p(w_jv_1)\\ \vdots\\ p(w_jv_{d^*})
   \end{pmatrix}
= W_ae_j = \begin{pmatrix} p(w_jav_1)\\ \vdots\\ p(w_jav_{d^*})
   \end{pmatrix}.
\end{split}
\end{equation}
\noindent $|v|\to|v|+1$: Let $\tilde{v}=av$ where $a\in\S$.
\begin{equation}
\begin{split}
T_{\tilde{v}}\begin{pmatrix} p(w_jv_1)\\ \vdots\\ p(w_jv_{d^*})
   \end{pmatrix}
&= T_vT_a\begin{pmatrix} p(w_jv_1)\\ \vdots\\ p(w_jv_{d^*})
   \end{pmatrix}
\stackrel{|v|=1}{=}T_v\begin{pmatrix} p(w_jav_1)\\ \vdots\\ p(w_jav_{d^*})
   \end{pmatrix}\\ 
&\stackrel{(*)}{=} \begin{pmatrix} p(w_jvav_1)\\ \vdots\\ p(w_jvav_{d^*})
   \end{pmatrix}
= \begin{pmatrix} p(w_j\tilde{v}v_1)\\ \vdots\\ p(w_j\tilde{v}v_{d^*})
   \end{pmatrix}
\end{split}
\end{equation}
where $(*)$ follows from the induction hypothesis.
\qed\\

\begin{Lem}
\label{l.generators.2}
For all $v,w\in\S^*$ such that $|w|,|v|\le\lceil\frac{n}{2}\rceil,|wv|\le n$
($T_v=T_{a_k}...T_{a_1},v=a_1...a_k\in\S^k$):
\begin{equation}
y^TT_v\begin{pmatrix} p(wv_1)\\ \vdots\\ p(wv_{d^*})
   \end{pmatrix}
= p(wv).
\end{equation}
\end{Lem}

{\bf Proof of lemma \ref{l.generators.2}}: Note that the columns in
$\cP_{p,\lfloor\frac{n}{2}\rfloor,\lceil\frac{n}{2}\rceil}$
resp.~$\cP_{p,\lceil\frac{n}{2}\rceil,\lfloor\frac{n}{2}\rfloor}$
referring to $w$ is contained in the span of the columns referring to
the $w_j$'s, according to the choice of the $w_j$.  Therefore, it
suffices to show the statement for $w=w_j$.  We do this by induction
on $|v|$:\par
\noindent $|v|=0$ ($v=\Box,T_{\Box}=Id$):
\begin{equation}
y^TT_{\Box}\begin{pmatrix} p(w_jv_1)\\ \vdots\\ p(w_jv_{d^*})
   \end{pmatrix} 
= y^T\begin{pmatrix} p(w_jv_1)\\ \vdots\\ p(w_jv_{d^*})
   \end{pmatrix} 
= p(w_j)
\end{equation}
follows from the choice of $y$.\par
\noindent $|v|\to|v|+1$: Let $\tilde{v}=av,a\in\S$.
\begin{equation}
\begin{split}
y^TT_{\tilde{v}}\begin{pmatrix} p(w_jv_1)\\ \vdots\\ p(w_jv_{d^*})
   \end{pmatrix}
&= y^TT_vT_a\begin{pmatrix} p(w_jv_1)\\ \vdots\\ p(w_jv_{d^*})
   \end{pmatrix}\\ 
& \stackrel{L.~\ref{l.generators.1}}{=}  
y^TT_v\begin{pmatrix} p(w_jav_1)\\ \vdots\\ p(w_jav_{d^*})
   \end{pmatrix} 
\stackrel{(*)}{=} p(wav) = p(w\tilde{v}) 
\end{split}
\end{equation}
where $(*)$ follows from the induction hypothesis.\qed\\

{\bf Proof of theorem \ref{t.generators} (cont.)}: Let $u\in\S^*$ such
that $|u|\le n$. Split $u=wv$ into two strings $w,v$ such that
$|w|,|v|\le \lceil\frac{n}{2}\rceil$. We compute
\begin{equation}
\begin{split}
y^TT_ux &= y^TT_vT_wx = y^TT_vT_w\begin{pmatrix} p(v_1)\\ \vdots\\ p(v_{d^*})
   \end{pmatrix} 
= y^TT_vT_wy^TT_vT_w\begin{pmatrix} p(\Box v_1)\\ \vdots\\ p(\Box v_{d^*})
   \end{pmatrix}\\
&\stackrel{L.~\ref{l.generators.1},|w\Box|\le\lceil\frac{n}{2}\rceil}{=} 
y^TT_v\begin{pmatrix} p(wv_1)\\ \vdots\\ p(wv_{d^*})
   \end{pmatrix}
\stackrel{L.~\ref{l.generators.2}}{=} p(wv) = p(u)
\end{split}
\end{equation} 
where we have replaced $v$ resp.~$w$ of lemma \ref{l.generators.1} by
$w$ resp.~$\Box$ here in order to obtain the fourth equation.\qed\\

Due to theorem \ref{t.generators}, invariants that are induced by
conditions $(a)$ and $(b)$ fully describe the finite-dimensional model
$\mathbf{g}_{n,d}$ for $n\ge 2d-1$, hence generate the ideal of model
invariants.

\begin{Xmp}
Consider 
\begin{equation*}
\cP_{p,4,2} = 
  \begin{pmatrix} 
       p(\Box) & p(0) & p(1)   & p(00)   & p(01)   & p(10)   & p(11)   \\
       p(0)  & p(00)  & p(10)  & p(000)  & p(010)  & p(100)  & p(110)  \\
       p(1)  & p(01)  & p(11)  & p(001)  & p(011)  & p(101)  & p(111)  \\
       p(00) & p(000) & p(100) & p(0000) & p(0100) & p(1000) & p(1100) \\
       p(01) & p(001) & p(101) & p(0001) & p(0101) & p(1001) & p(1101) \\
       p(10) & p(010) & p(110) & p(0010) & p(0110) & p(1010) & p(1110) \\
       p(11) & p(011) & p(111) & p(0011) & p(0111) & p(1011) & p(1111) \\
  \end{pmatrix}
\end{equation*}
where $\S=\{0,1\}$. Condition $(a)$ then translates to the only equation
\begin{equation*}
\det
  \begin{pmatrix} 
       p(\Box) & p(0) & p(1)\\   
       p(0)  & p(00)  & p(10)\\  
       p(1)  & p(01)  & p(11)  
  \end{pmatrix} = 0.
\end{equation*}
The column conditions in $(b)$ can be stated as follows:
\begin{multline*}
\begin{pmatrix} 
        p(00)   \\
        p(000)  \\
        p(001)  \\
        p(0000) \\
        p(0001) \\
        p(0010) \\
        p(0011) \\
  \end{pmatrix}, 
\begin{pmatrix} 
        p(01)   \\
        p(010)  \\
        p(011)  \\
        p(0100) \\
        p(0101) \\
        p(0110) \\
        p(0111) \\
  \end{pmatrix},
\begin{pmatrix} 
       p(10)   \\
       p(100)  \\
       p(101)  \\
       p(1000) \\
       p(1001) \\
       p(1010) \\
       p(1011) \\
  \end{pmatrix},
\begin{pmatrix} 
       p(11)   \\
       p(110)  \\
       p(111)  \\
       p(1100) \\
       p(1101) \\
       p(1110) \\
       p(1111) \\
  \end{pmatrix}\\
\in\spann\{
\begin{pmatrix} 
       p(\Box)\\
       p(0)  \\
       p(1)  \\
       p(00) \\
       p(01) \\
       p(10) \\
       p(11)
  \end{pmatrix},
\begin{pmatrix} 
       p(0)   \\
       p(00)  \\
       p(01)  \\
       p(000) \\
       p(001) \\
       p(010) \\
       p(011) 
  \end{pmatrix}
\begin{pmatrix} 
       p(1)   \\
       p(10)  \\
       p(11)  \\
       p(100) \\
       p(101) \\
       p(110) \\
       p(111) \\
  \end{pmatrix}\}
\end{multline*}
The row conditions are completely analogous to the column conditions.
\end{Xmp}

\begin{sloppypar}
Clearly, invariants induced by $(b)$
refer to polynomial rings
\begin{equation}
K[X_{ij},Y_i,1\le i\le M, 1\le j\le N]
\end{equation}
and the smallest varieties therein that contain all points
$x_{ij},y_i$ such that $(y_1,...,y_M)$ is linearly dependent of
$(x_{11},...,x_{M1}),...,(x_{1N},...,x_{MN})$. The Zariski closure of
the image of $\mathbf{g}_{n,d}$ being an irreducible variety leads us
to the following conjecture.\end{sloppypar}

\begin{Conj}
Let $n\ge 2d-1$.
\begin{equation*}
(p(v))_{v\in\S^n}\in\overline{\image(\mathbf{g}_{n,d})}
\end{equation*}
if and only if
\begin{equation*}
\det[p(w_jv_i)]_{1\le i,j\le d+1} = 0
\end{equation*}
for all choices of words $v_1,...v_{d+1},w_1,...w_{d+1}$ 
such that $|w_jv_i|\le n$.
\end{Conj}

\begin{Rem}
\label{rem.ooms}
The finite-dimensional models have to be handled with certain {\em
  care}.  Even if a (unconstrained) probability distribution is in the
image of $\mathbf{g}_{n,d}$ the finite-dimensional string function
giving rise to it might not be an (unconstrained) stochastic process,
meaning that the string function is not necessarily non-negative,
since values referring to longer strings as computed according to
(\ref{eq.oom}) might be negative.  It is one of the {\em big open
  problems} of the theory of finite-dimensional processes how to
algorithmically determine whether a set of matrices as in
(\ref{eq.oom}) gives rise to a non-negative string function.
\end{Rem}

\section{String Length Complexity}
\label{sec.finitedet}

In this section, we will prove a set-theoretical theorem an
ideal-theoretical counterpart of which would yield, as a corollary, a
proof of conjecture 11.9, \cite{Bray05}. The theorem may be of
interest in its own right, as the assumptions to be met by the models
under considerations are fairly mild.\par Roughly speaking, an
ideal-theoretical extension of the theorem would be about how to lift
sets of generators for models describing distributions over strings of
length $n$ to generators for distributions over strings of length
$n+1$, given that $n$ is greater than the {\em string length
  complexity} of the underlying models.

In the following,  
\begin{equation}
\cM\subset\R^{\S^*}
\end{equation}
is a class of USSFs. 

\begin{Def}
Let $\cM\subset\R^{\S^*}$ be a class of USSFs.  We define the {\em
  string length complexity} of $\cM$ to be
\begin{multline}
\SLC(\cM):=\\
\inf\{N\in\N\;|\;
p_1,p_2\in\cM:\;(p_1)_{|\S^n}= (p_2)_{|\S^n}\;\Rightarrow\;p_1 = p_2\}.
\end{multline}
That is, members of $\cM$ are uniquely determined by their
distributions over strings of length $\SLC(\cM)$.
\end{Def}

Given a class of USSFs, let
\begin{equation}
\cM_n :=\{(p(v))_{v\in\S^n}\,|\,p\in\cM\}
\end{equation}
be the set of distributions over strings of length $n$ that are induced
by the members of $\cM$. In case of $\SLC(\cM)=n$ the map
\begin{equation}
\begin{array}{rccc}
\pi_{\S^n}: & \cM & \longrightarrow & \cM_n\\
           & p   & \mapsto         & p_{|\S^n} = (p(v))_{v\in\S^n}
\end{array}
\end{equation}
is one-to-one.

\begin{Thm}
\label{t.degree}
Let $\cM$ be a class of unconstrained random processes such that
\begin{itemize}
\item[(i)]
  \begin{equation}
    \label{eq.finiteslc}
    \SLC(\cM) \le n-1 < \infty.
  \end{equation}
\item[(ii)]
  \begin{equation}
    \label{eq.ppa}
    p\in\cM\quad\Rightarrow\quad \forall a\in\S:\; p^a\in\cM.
  \end{equation}
\end{itemize}

Then it holds that
\begin{equation}
\label{eq.lift}
(p(u),u\in\S^{n+1})\in\cM_{n+1}\quad\Leftrightarrow\quad
\begin{cases}
(p(av),v\in\S^n)\in\cM_n & \forall a\in\S\\
(p(v),v\in\S^n)\in\cM_n & 
\end{cases}
\end{equation}
where $p(v)=\sum_{a\in\S}p(va)$.

\end{Thm}

\begin{Rem}
Theorem \ref{t.degree} is meant to be a first step to obtain an
analogous theorem resulting from replacing $\cM_n,\cM_{n+1}$ by their
Zariski closures $\overline{\cM_n},\overline{\cM_{n+1}}$.  Generators
for the ideal of invariants of $\overline{\cM_{n+1}}$, given
generators for the ideal of invariants of $\overline{\cM_n}$, could be
obtained by the following idea. If $h\in\C[X_v,v\in\S^n]$ is one of the
generators for $\overline{\cM_n}$ where the $X_v$ are indeterminates
for the probabilities $p(v),v\in\S^n$, one obtains $|\S|+1$ generators
for $\overline{\cM_{n+1}}$ by replacing the indeterminates
$X_v,v\in\S^n$ by indeterminates $X_{av},v\in\S^n$ for all $a\in\S$
which results in new generators
\begin{equation}
h_a\in\C[X_{av},v\in\S^n]\subset\C[X_u,u\in\S^{n+1}]
\end{equation}
as well as replacing each $X_v$ by the polynomials
$\sum_aX_{va}\in\C[X_u,u\in\S^{n+1}]$ resulting in another generator
\begin{equation}
h_+\in\C[\sum_aX_{va},v\in\S^n]\subset\C[X_u,u\in\S^{n+1}]. 
\end{equation}
The theorem would state that the generators obtained by this procedure
generate the ideal of invariants of $\cM_{n+1}$.\par
Note that in particular the maximum degree of the generators of $\overline{\cM_{n+1}}$
would be at most that of $\overline{\cM_n}$.
\end{Rem}

\Pf ``$\Rightarrow$'': From (\ref{eq.ppa}) we obtain that
$(p^a(v),v\in\S^n)\in\cM_n$ for each $a\in\S$. 
The second part is just the trivial observation that
$(p(u),u\in\S^{n+1})\in\cM_{n+1}$ implies $(p(v),v\in\S^n)\in\cM_n$.\\

\noindent ``$\Leftarrow$'': 
From the second case in (\ref{eq.lift}) we obtain that 
$(p(v),v\in\S^n)\in\cM_n$. As elements of $\cM$ are uniquely determined
by their values for strings of length at least $m$ and $n\ge m+1$ we
obtain a USSF $\tilde{p}\in\cM$ such that
\begin{equation}
\label{eq.ptildep}
p(v) = \tilde{p}(v) \quad\text{ for all }v\in\S^n.
\end{equation}
It remains to show that also
\begin{equation}
p(w) = \tilde{p}(w) \quad\text{ for all }w\in\S^{n+1}
\end{equation}
which amounts to showing that
\begin{equation}
\label{eq.pavtildepav}
p(av) = \tilde{p}(av)=\tilde{p}^a(v) \quad\text{ for all }(a,v)\in\S\times\S^{n}.
\end{equation}
We further observe that  
\begin{equation}
(p^a(v),v\in\S^n)\in\cM_n
\end{equation}
for all $a\in\S$, because of $n\ge m+1> m$, implies the existence of a
unique $q^a\in\cM$ s.t.~
\begin{equation}
\label{eq.qapl}
q^a(v) = p^a(v)\quad\text{ for all }v, |v|\le n.
\end{equation}
As $\tilde{p}\in\cM$, we have that $\tilde{p}^a\in\cM$ for all $a\in\S$,
due to (\ref{eq.ppa}). Moreover, for $u\in\S^{n-1}$,
\begin{equation}
\tilde{p}^a(u) = \tilde{p}(au) \stackrel{(\ref{eq.ptildep})}{=} p(au)
\stackrel{(\ref{eq.qapl})}{=} q^a(u).
\end{equation}
As $n-1\ge m$ and $\tilde{p}^a\cM$ and $q^a\cM$ coincide on strings of
length $n-1\ge m$, we obtain
\begin{equation}
\label{eq.qatildep}
\tilde{p}^a = q^a
\end{equation}
because of $(i)$. We finally compute
\begin{equation}
p(av) = p^a(v)\stackrel{(\ref{eq.qapl})}{=} q^a(v) \stackrel{(\ref{eq.qatildep})}{=}
\tilde{p}^a(v) = \tilde{p}(av)
\end{equation}
which establishes (\ref{eq.pavtildepav}).\qed\\

\subsection{Finite-Dimensional Models}

Theorem \ref{t.degree} applies for the finite-dimensional
models. Eq.~\ref{eq.finiteslc} is established by theorem \ref{t.unique} in 
subsection \ref{ssec.dim} whose statement is that finite-dimensional
processes $p$ of dimension at most $d$ are uniquely determined by the
values $p(v),|v|=2d-1$.

In terms of the language introduced here, we can restate theorem \ref{t.unique}
as follows.

\begin{Thm}
\label{t.dimslc}
Let 
\begin{equation*}
\cM_d:=\{p\in\R^{\S^*};|\;p\text{ is USSF and } \dim p\le d\}
\end{equation*}
be the class of unconstrained processes of dimension at most $d$.
Then it holds that
\begin{equation*}
\SLC(\cM_d) = 2d-1.
\end{equation*}
\end{Thm}

Furthermore observe that 
\begin{equation}
(p^a)^w(v) = p^a(wv) = p(awv) = p^{aw}(v)
\end{equation}
for all $a\in\S,v,w\in\S^*$ 
which translates to 
\begin{equation}
(p^a)^w = p^{aw}.
\end{equation}
Hence the column space of $\cP_{p^a}$ is contained in that of $\cP_p$ which yields
\begin{equation}
\dim p^a\le\dim p
\end{equation}
as $\dim p$ is just the dimension of the column space of $\cP_p$.

This observation in combination with theorem \ref{t.dimslc} make the assumptions
of theorem \ref{t.degree} hold for $\cM_d$, which yields the following corollary.

\begin{Cor}
\label{c.uniquefinite}
Let $n\ge 2d$. Then it holds that
\begin{equation}
\label{eq.finitelift}
(p(u),u\in\S^{n+1})\in\image\mathbf{g}_{n+1,d}\quad\Leftrightarrow\quad
\begin{cases}
(p(av),v\in\S^n)\in\image\mathbf{g}_{n,d} & \forall a\in\S\\
(p(v),v\in\S^n)\in\image\mathbf{g}_{n,d} & 
\end{cases}
\end{equation}
 \end{Cor}

Again, an analogous ideal-theoretical result referring to the Zariski
closures of $\image\mathbf{g}_{n,d},\image\mathbf{g}_{n+1,d}$ would
yield that the maximum degree of the generators would not increase for
$n\ge 2d$.

\section{Hidden Markov Models}
\label{sec.hmms}

In the following, let
\begin{equation}
\begin{array}{rccc}
\mathbf{f}_{n,l}: & \C^{l(l-1)+l(|\S|-1)+l} & \longrightarrow & \C^{|\S|^n}\\
                 &((T_a=A^TO_a)_{a\in\S}),x)&\mapsto&(\tr T_{a_n}...T_{a_1}x(1,...,1)^T)_{v=a_1...a_n\in\S^n}.
\end{array}
\end{equation}
where $A$ and the $O_a$ as in example \ref{x.hmm}, be the polynomial
map associated with the unconstrained (constrained if and only
if $\sum_{i=1}^lx_i=1$) hidden Markov model referring to hidden Markov
models acting on $l$ hidden states and distributions over strings of
length $n$, as described in \cite{Pachter05}. \par The following
theorem of Heller resulted from the attempts set off in the late 50's
\cite{Gilbert59,Dharma63a,Dharma63b,Dharma65} to give novel
characterizations of hidden Markov processes. Many of those results
are based on the idea that HMPs have finite dimension, which was
noticed earlier in that series of papers without explicitly stating
it. We give a version of Heller's theorem that is adapted to the
language in use here. Heller's version is formulated in the language
of homological algebra---without string functions and Hankel
matrices. In his paper, discrete random processes are viewed as
modules over certain rings. This language later has never been used in
the theory of stochastic processes or related areas, probably as the
required amount of prior knowledge unfamiliar to statisticians and
probabilists is high. In the following we define
\begin{equation}
C_p:=\spann\{p^w\,|\,w\in\S^*\}
\end{equation}
to be the column space of the Hankel matrix $\cP_p$ of a string
function $p$.

\begin{Thm}[Heller, 1965]
\label{t.heller}
A string function $p:\S^*\to\R$ is associated with a (unconstrained) hidden Markov
process if and only if there are (U)SSFs
$p_i\in C_p,i=1,...,l$ s.t.~
\begin{enumerate}
\item[(a)] $p\in\cone\{p_i\,|\,i=1,...,l\}$,
\item[(b)] $\forall w\in\S^*:\;(p_i)^w\in\cone\{p_i\,|\,i=1,...,l\}$.
\end{enumerate}
\end{Thm}

Note first that this again points out that hidden Markov processes $p$
are finite-dimensional as $C_p\subset\spann\{p_i\,|\,
i=1,...,l\}$ hence $\dim p\le l$.  Note further that $(a)$ in
combination with $(b)$ implies that $p^v\in\cone\{p_i\}$ for all
$v\in\S^*$ which renders $\cone\{p_i\}$ to be full-dimensional. It is closed
due to being polyhedral and pointed due to being generated by SSFs
which are strictly positive string functions. Collecting properties
results in $\cone\{p_i\}$ being a proper, polyhedral cone.\\

Given a hidden Markov process $p$, the $p_i$ can be obtained as the
random processes starting from the hidden states (i.e.~having initial
probability distribution $e_i$). The other direction requires more
work. A translation of Heller's proof \cite{Heller65} to the language
of string functions can be found in \cite{Schoenhuth06}. A rather
straightforward consequence of Heller's theorem is the following
corollary.

\begin{Cor}
Let $p$ be a USSF of dimension of at most $2$. Then $p$ is associated
with an unconstrained hidden Markov process acting on $2$ hidden states.
\end{Cor}

{\em Proof Sketch}: As all $p^w\ge 0$ the cone generated by all column 
vectors
\begin{equation}
\cone\{p^w\,|\,w\in\S^*\}
\end{equation}
is pointed hence its closure is generated by its extremal rays. In two
dimensions this is equivalent to the closure of $\cone\{p^w\,|\,v\in\S^*\}$
being polyhedral. It's a routine exercise to check for the assumptions
of Heller's theorem to hold for this cone.\qed\\

One might be tempted to infer that the ideal of model invariants of
$\mathbf{f}_{n,2}$ can be computed by computing the invariants of the
$2$-dimensional model, as provided by theorem
\ref{t.generators}. However, a $2$-dimensional process need not be
associated with a hidden Markov process acting on $2$ hidden states.
According to the proof of theorem \ref{t.heller}, one might need up
to $2|\S|$ many hidden states to describe an arbitrary $2$-dimensional
process by means of a hidden Markov parameterization.\\

\subsection{Degree of Invariants}

Heller's theorem gives rise to an application of theorem
\ref{t.degree} to hidden Markov processes where $n\ge 2l$. Assumption
$(i)$ of theorem \ref{t.degree} is met since hidden Markov processes
on $l$ hidden states, as finite-dimensional random processes of
dimension $\le l$, are determined by their distributions over the
strings of length $2l-1$. Assumption $(ii)$ is met due to Heller's
theorem.\footnote{Proofs for this can also be formulated in terms of
  the hidden Markov processes' parameterizations. However, such proofs
  are lengthy and technical exercises.}  The only thing one has to be
aware of is that the dimension of the column space of $\cP_{p^a}$ can
be lower than that of $\cP_p$ itself. In this case, one obtains the
necessary cone generators by projecting $C_p$ onto $C_{p^a}$ (we
recall that $C_{p^a}\subset C_p$).  In sum, the class of unconstrained
hidden Markov processes meet the assumptions of theorem
\ref{t.degree}, which yields

\begin{Cor}
\label{c.uniquehmm}
Let $n\ge 2d$. Then it holds that
\begin{equation}
  \label{eq.hmmlift}
  (p(u),u\in\S^{n+1})\in\image\mathbf{f}_{n+1,l}\quad\Leftrightarrow\quad
  \begin{cases}
    (p(av),v\in\S^n)\in\image\mathbf{f}_{n,l} & \forall a\in\S\\
    (p(v),v\in\S^n)\in\image\mathbf{f}_{n,l} & 
  \end{cases}
\end{equation}
\end{Cor}

Note that an ideal-theoretic equivalent of theorem \ref{c.uniquehmm}
would yield a proof of conjecture 11.9 from \cite{Bray05} as a corollary.
However, an ideal-theoretical equivalent of theorem \ref{c.uniquehmm} would 
be a stronger result:

\begin{Conj}
Let $\mathbf{f}_{n,l}$ be the unconstrained hidden Markov model for
$l$ hidden states and strings of length $n$. Then the maximum degree
of the invariants $d(n,l)$ of $\mathbf{f}_{n,l}$ does not increase for
$n\ge 2l$, that is,
\begin{equation}
\quad ...d(n+1,l) \le d(n,l)\le d(n-1,l)\le ...\le d(2l,l).
\end{equation}
\end{Conj}

As $d(5,2)=1$ (see \cite{Bray05}, table 11.1 (?)), we would obtain
that $d(n,l)=1$ for $n\ge 5$, that is, the ideal of invariants would
be generated by linear equations exclusively.

\section{The Markov model}
\label{sec.markov}

In the following, let (U)SSFs $p$ be induced by Markov chains. That is,
\begin{equation}
p(v=a_1...a_n) = \pi(a_1)\prod_{i=2}^nM_{a_{i-1}a_i}
\end{equation}
where $\pi\in\R^{\S}$ is a strictly positive vector (with entries
not necessarily summing up to one in case of a USSF $p$) and
$M\in\R^{\S^2}$ is a matrix with the entries of a row summing up to
one. Moreover, in this section, let
\begin{equation}
\begin{array}{rccc}
\mathbf{f}_{n,l=|\S|}: & \C^{l+l(l-1)} & \longrightarrow & \C^{|\S|^n}\\
                      & (\pi,M) & \mapsto & (\pi(a_1)\prod_{i=2}^nM_{a_{i-1}a_i})_{v=a_1...a_n\in\S^n}.
\end{array}
\end{equation}
be the polynomial map (associated with the Markov model in case of
$\pi,M$ being in accordance with the laws from above) with
alphabet $\S$ on strings of length $n$. In the language of string
functions and Hankel matrices, we have the following theorem.

\begin{Thm}
\label{t.markov}
A (U)SSF $p$ is associated with a Markov chain iff
\begin{equation}
\label{eq.markov}
\forall a\in\S:\quad\dim\spann\{p^{va}\,|\,v\in\S^*\}\le 1.
\end{equation}
\end{Thm}

A proof can be found in \cite{Schoenhuth06}, for example.\par
This can be straightforwardly exploited to obtain invariants of
$\mathbf{f}_{n,l}$. 

\begin{Thm}
\label{t.markovinvar}
Let $(p(v),v\in\S^n)$ be a (unconstrained) probability distribution
such that $n\ge 2|\S|-1$.  Then $(p(v),v\in\S^n)$ lies in the image of
$\mathbf{f}_{n,l=|\S|}$ if and only if
\begin{equation}
\label{eq.markovhankel}
\det\begin{bmatrix} p(vau) & p(wau)\\
                    p(vau') & p(wau')
    \end{bmatrix} = 0
\end{equation}
for all choices $u,u',v,w\in\S^*,a\in\S$ such that
$|vau|,|vau'|,|wau|,|wau'|\le n$ and, as usual,
$p(v):=\sum_{w\in\S^{n-|v|}}p(vw)$ for strings $v$ such that $|v|<n$.
\end{Thm}

\Pf ``$\Rightarrow$'' is obvious as for a Markov chain $p$,
(\ref{eq.markovhankel}) is a necessary consequence of
(\ref{eq.markov}) in theorem \ref{t.markov}.\par ``$\Leftarrow$''
Clearly, (\ref{eq.markovhankel}) implies the assumptions
(\ref{eq.genaa}) and (\ref{eq.condb}) of theorem
\ref{t.generators} to hold, which yields that $(p(v),v\in\S^n)$ lies
in the image of the finite-dimensional model. We thus find, by means
of algorithm \ref{a.ooms}, $(T_a)_{a\in\S},x,y$ such that the
probabilities $p(v)$ for all $v$ up to length $n\ge 2|\S|$ can be
computed according to (\ref{eq.oom}).  Note that $T_a$ maps $p^v$ onto
$p^{va}$ where $p^v,p^{va}$ are identified with a coordinate
representation induced by the basis of the column spaces that one has
found according to algorithm \ref{a.ooms} (see remark \ref{rem.ooms}).
In this sense, (\ref{eq.markovhankel}) translates to 
\begin{equation}
\dim\image T_a\le 1
\end{equation}
for all $a\in\S$. Clearly, this implies (\ref{eq.markov}) of theorem
\ref{t.markov} from which the assertion follows.\qed\\

\begin{Rem}
While the assumption $n\ge 2|\S|-1$ helps to give a rather concise proof of
theorem \ref{t.markovinvar}, we feel that it is not a necessary
requirement. However, inference of Markov chain parameters giving rise
to probability distributions $(p(v),v\in\S^n)$ for which the
determinantal invariants (\ref{eq.markovhankel}) apply is a much more
technical undertaking. Moreover, it seems that some (potentially more
involved) pecularities have to be resolved.
\end{Rem}

\section{Trace algebras}
\label{sec.tracealg}

In this section, we will draw some connections between trace algebras
and the theory of finite-dimensional string functions.  For a rigorous
introduction to trace algebras see \cite{Drensky07}. We recall that in
Bernd's preprint \cite{Bernd08} the quartic hidden Markov model
invariant listed in \cite{Bray05} could be identified as a relation
between trace polynomials.\par Here, we shall try to shed some light
on the general relationships between trace algebras and
finite-dimensional models. In terms of the language of trace
algebras, we will derive some defining relations for the trace
algebras.\par Therefore, we introduce the following definition.

\begin{Def}
A string function $p:\S^*\to\R$ is called {\em traceable of order $r$}
if there are matrices $X_a\in\R^{r\times r},a\in\S$ such that
\begin{equation}
p(v=a_1...a_n) = \tr X_{a_n}...X_{a_1}.
\end{equation}
\end{Def}

Traceable string functions are finite-dimensional, as can be seen
by application of a simple lemma.

\begin{Lem}
\label{l.sumdim}
Let $p_i,i=1,...,k$ be string functions of dimensions $d_i$. Let
$p:=\sum_{i=1}^kp_i$.  Then it holds that
\begin{equation}
\dim p \le \sum_{i=1}^kd_i.
\end{equation}
\end{Lem}

This gives rise to

\begin{Thm}
\label{t.tracedim}
Let $p\in\R^{\S^*}$ be traceable of order $r$.  Then
\begin{equation}
\dim p\le r^2.
\end{equation}
\end{Thm}

\Pf Let $p_i\in\R^{\S^*},i=1,...,r$ be defined by 
\begin{equation}
p_i(v=a_1...a_n) := \tr X_{a_n}...X_{a_1}e_ie_i^T.
\end{equation}
From theorem \ref{t.oom} we obtain $\dim p_i \le r$.
As $Id = \sum_ie_ie_i^T$, which yields 
\begin{equation}
p = \sum p_i
\end{equation}
the assertion follows from application of lemma \ref{l.sumdim}.\qed\\

If the identity matrix $Id=\sum_ie_ie_i^T$ was presentable in the form
$Id=xy^T$ itself, traceable string functions of order $r$ would be of
dimension at most $r$, as given by theorem \ref{t.oom}. As this is
not the case, there are traceable string functions of order $r$ whose
dimension is larger than $r$. Moreover, not every string function of
dimension $r^2$ seems to be traceable. However, an example of that
kind is yet to be delivered.\\

The consequences of theorem \ref{t.tracedim} for the theory of trace
algebras are that invariants which can be computed for the
$r^2$-dimensional models $\mathbf{f}_{n,r^2}$ also apply as defining
relations for the trace algebras generated by all trace polynomials
\begin{equation}
\tr (X_{i_n}...X_{i_1}), 1\le i_j\le d, n\ge 0.
\end{equation}

The exact relationships between trace algebras, hidden Markov as well
as the finite-dimensional models are yet to be determined.

\section{Open Questions}
\label{sec.openq}

\begin{enumerate}
\item Theorem \ref{t.heller} characterizes hidden Markov chains within
  the theory of finite-dimensional random processes. Determine
  invariants that correspond to this characterization.
\item Determine the relationships between trace algebras and the
  models under consideration here in more detail.
\item Deliver a proof for a more general version of theorem
  \ref{t.markovinvar}, as discussed above.
\item Determine the peculiarities of differences between the
  two-dimensional models and the hidden Markov models for $2$ hidden
  states.
\item Tropicalization of Teichm\"uller spaces (see \cite{Bernd08})?
\end{enumerate}

\end{document}